\documentclass[12pt,a4paper,reqno,dvipsnames]{amsart}

\usepackage[utf8]{inputenc}
\usepackage[english]{babel}
\usepackage[normalem]{ulem} %para usar o comando "\sout" (strike out) que permite "riscar" sentenças

\usepackage{enumerate}
\usepackage{xypic}
\usepackage{amsmath}
\usepackage{amsthm}
\usepackage{amsfonts}
\usepackage{amssymb}
\usepackage{graphicx,color}
\usepackage{geometry}
\usepackage{wrapfig}
\usepackage{setspace}
\usepackage{multicol}
\usepackage{xcolor}
\usepackage{soul} %% para usar com o comando \hl{...} que destaca texto
\usepackage{hyperref}

\newtheorem{theorem}{Theorem}
\newtheorem{lemma}[theorem]{Lemma}
\newtheorem{proposition}[theorem]{Proposition}
\newtheorem{definition}[theorem]{Definition}
\newtheorem{remark}[theorem]{Remark}
\newtheorem{corollary}[theorem]{Corollary}
\newtheorem{example}[theorem]{Example}

\newtheorem*{theorem*}{Theorem}
\newtheorem*{proposition*}{Proposition}

\newcommand{\Hom}{\operatorname{Hom}}

\newcommand{\K}{\Bbbk}

\newcommand{\mC}{\mathcal{C}}



\title[Partial Hopf coactions on nonunital algebras]{PARTIAL COACTIONS OF HOPF ALGEBRAS ON NONUNITAL ALGEBRAS AND CATEGORIES}
\author[M. \ M. \ Alves]{Marcelo Muniz Alves}
\address{Departamento de Matem\'atica, Centro Polit\'ecnico, Universidade Federal do Paran\'a, Curitiba - PR, Brasil}
\email{marcelomsa@ufpr.br}

\author[T. L. Ferazza]{Tiago Luiz Ferrazza}
\address{Universidade Estadual do Paran\'a  - Campus de Paranagu\'a, Rua Comendador Correa Junior, 117, Paranaguá - PR, Brasil}
\email{tiago.ferrazza@unespar.edu.br}

\thanks{This work was partially financed  by  Coordena\c{c}\~ao de Aperfeiçoamento de Pessoal de N\'ivel Superior – Brasil (CAPES) – Finance Code 001, Fundação de Amparo \`a Pesquisa do Estado de Minas Gerais (FAPEMIG) under Universal Grant APQ-03491-25, and National Council for Scientific and Technological Development – CNPq under Universal Grant 405540/2023-0.}

\begin{document}

\begin{abstract}
In this paper, we introduce and investigate partial coactions of a Hopf algebra on an associative, possibly nonunital algebra and on $\K$-linear categories. We prove that every partial coaction has a globalization in a broad sense. We study in detail the case of a partial coaction on an algebra with local units; we prove that such a coaction is always induced from a global coaction, and we show that there exists a naturally associated $\K$-linear category endowed with a partial coaction.
\end{abstract}

\maketitle

\section{Introduction}

Partial actions of groups on algebras were introduced in \cite{Circle.actions}, in the field of operator algebras, as a means of expressing a $C^*$-algebra with a circle action as a generalized skew group algebra over the integers. This construction was later generalized to any group and any associative algebra (see, for instance, \cite{dok}).

In the literature, partial actions of groups are usually defined on associative nonunital algebras, and a great part of the recent developments around partial actions and graded algebras still happens in the context of nonunital algebras \cite{abadie2024strong,dok3,dok}. Partial actions and coactions of Hopf algebras followed another path and were first defined for unital algebras; they first appear in \cite{caenepeel}, motivated by the Galois theory of partial group actions \cite{ferrero}.

A main question in partial (co)actions is the existence of a globalization: is every partial (co)action isomorphic to one obtained from a global (co)action? If not, can one characterize those that are globalizable? A classical result in this direction appears in \cite{dok}, where it is shown that a partial group action on an algebra is globalizable if and only if it is a unital action. On the other hand, every partial Hopf (co)action on a unital algebra is globalizable \cite{Alves1,alves3}. In this paper, we show that partial Hopf coactions on nonunital algebras are always globalizable, although via a weaker form of globalization. This result is strengthened in the case of algebras with local units.

Our motivation to extend partial actions and coactions to nonunital algebras comes from the related theory of actions and coactions on $\K$-linear categories, also known as $H$-module categories and $H$-comodule categories, respectively.

$H$-module categories were introduced in \cite{cibils-solotar-galois} as a means of studying problems arising from representation theory of associative algebras: they provide a means of unifying constructions of Galois coverings for $\K$-linear categories via a smash product of a Hopf algebra by a $\K$-linear category.

$H$-comodule categories were introduced in \cite{herscovich2007hochschild}, modeled after $G$-gradings of categories, where the authors introduce a notion of Hopf--Galois extension for linear categories. The authors obtain a spectral sequence relating the Hochschild--Mitchell cohomology of the coinvariant subcategory to that of the original category, extending results of Lorenz \cite{lorenz} on strongly graded algebras and Stefan \cite{stefan1995} on Hopf--Galois extensions of algebras.

These Hopf--Galois extensions of categories were further studied in \cite{anca2,anca1}, where a theory of cleft extensions of categories was developed. In \cite{balodi2020cohomology}, the authors extend the results on the cohomology of relative Hopf modules of \cite{caenepeel2005cohomology} to Hopf comodule categories.

These topics are related to actions and coactions on algebras with local units via the ``matrix algebra construction'': given a $\K$-linear category $\mathcal{C}$, the direct sum of Hom spaces $a(\mathcal{C}) = \oplus_{x,y} {}_y\mathcal{C}_x$ has a canonical algebra structure, a canonical system of local units, and $H$-actions (coactions) on $\mathcal{C}$ induce $H$-actions (coactions) on $a(\mathcal{C})$ that preserve the system of local units.

Partial actions on nonunital algebras were investigated in \cite{ferrazza1}. The authors defined a globalization of such a partial action and proved that there always exists a globalization. With respect to $\K$-linear categories, it was shown in \cite{alves2} that a partial $H$-action on a $\K$-linear category $\mathcal{C}$ gives rise to a partial action on the algebra $a(\mathcal{C})$. Conversely, it is shown that a class of partial actions on algebras with local units corresponds to a partial action on an associated category.

Partial (co)actions of $C^*$-bialgebras were introduced in \cite{timmermann}, where the authors investigate the globalization problem. Partial coactions of multiplier Hopf algebras on nonunital algebras were introduced and investigated in \cite{azevedo2021partial}, where the authors link partial coactions to a Morita context associated to a partial action. The partial coactions considered in the present work are not particular cases of these previous ones, since we do not introduce idempotent multipliers in order to define them.

In this paper, we introduce and investigate partial coactions of a Hopf algebra on an associative (possibly nonunital) algebra and on $\K$-linear categories. We begin by introducing partial coactions on associative algebras in Section \ref{section2;partialcoactions}. We show that this definition comes essentially from considerations on a bimodule structure of $A\otimes H$, exactly as in the first paper on partial actions and coactions by Caenepeel and Janssen \cite{caenepeel}.

In Section \ref{section3:globalization} we introduce the concept of globalization for a partial coaction and we prove the main result of this paper: every partial coaction is isomorphic to an ideal of an $H$-comodule algebra which is generated by this ideal as an $H$-comodule (Theorem \ref{theorem.globalization}). This $H$-comodule algebra is the globalization of the original partial $H$-comodule algebra. 
\begin{theorem*}
Let $\rho: A\to A\otimes H$ be a symmetrical partial coaction, and let $B$ be the subcomodule algebra of $A\otimes H$ generated by $\rho(A)$. Then the pair $(B,\rho)$ is a globalization for $\rho$.
\end{theorem*}
We show that if $H$ is a finite-dimensional Hopf algebra then the minimal globalization of a partial action coincides with a globalization of the induced partial coaction (Proposition \ref{prop2}). We apply this result to partial good gradings. 

In Section \ref{section4.local.units} we investigate partial coactions on algebras with local units. We begin by considering the following scenario: $B$ is an $H$-comodule algebra (via $\gamma: B \to B \otimes H$) and $A$ is an ideal of $B$ which, as an algebra, has a system of local units. We show that, under a mild compatibility hypothesis, it is possible to induce a partial coaction on $A$ (Proposition \ref{coacparcialideal}). 
Roughly, we define a partial coaction $\rho: A \to A \otimes H$ pointwise with the aid of the system of local units: for each $a \in A$, we pick up an appropriate local unit $e \in A$ and we define $\rho(a)$ as 
\[
\rho(a) = (e \otimes 1_H)\gamma(a).
\]
The compatibility condition is used to patch together these expressions into a well-defined mapping partial coaction on $A$.

We then show that if $\rho: A \to A \otimes H$ is a partial coaction and $A$ has local units, then the globalization obtained in Theorem \ref{theorem.globalization} satisfies the compatibility condition. Hence, for algebras with local units, not only is every partial coaction globalizable, but every partial coaction is \emph{induced} by a global coaction (Proposition \ref{prop:glob.local.units} and Theorem \ref{thm:global.local.units}).

In the same section, we introduce partial coactions on $\K$-categories. We also introduce a class of coactions on algebras with a distinguished system of local units $S$, the $S$-categorizable partial coactions. We show how to obtain an $S$-categorizable coaction from a partial $H$-comodule category, and conversely. We conclude with a result that establishes a correspondence at the level of globalizations (Theorems \ref{thm:global.S-categorical.to.category} and \ref{thm:global.category.to.S-categorical.}).

%%%%%%%%%%%

\section{Partial coactions on nonunital algebras \label{section2;partialcoactions}}

In this section we introduce the concept of a partial coaction of a Hopf algebra $H$ on an associative, possibly nonunital, algebra $A$, extending the notion of partial coaction on unital algebras given in \cite{caenepeel}. As in \cite{caenepeel}, the definition of partial coaction that we propose is related to the $A$-bimodule structure of $A\otimes H$ and to the $A$-linearity of the map $I\otimes \varepsilon: A\otimes H\to A$.

First, let us recall the definition of partial coactions of Hopf algebras on unital algebras.

\begin{definition}[\cite{caenepeel}] \label{def:partial.coaction.unital} Let $H$ be a Hopf algebra and $A$ be a unital algebra. A linear map ${\rho : A \longrightarrow A \otimes H}$ is a partial $H$-coaction on $A$ if the following holds:
\begin{enumerate}[\normalfont(i)]
	\item $(I\otimes \varepsilon_H)\rho(a)=a$;
	\item $\rho(ab)=\rho(a)\rho(b)$;
	\item $(\rho \otimes I)\rho(a)=(\rho(1_A)\otimes 1_H)(I \otimes \Delta)\rho(a)$,
\end{enumerate}
for every $a,b \in A$. In this case, $A$ is called a \textit{partial $H$-comodule algebra}. If, additionally, $(\rho \otimes I)\rho(a)=[(I \otimes \Delta_H)\rho(a)](\rho(1_A)\otimes 1_H)$ for every $a\in A$, then $\rho$ is a symmetrical partial coaction. 
\end{definition}

Of course, if we wish to work with a nonunital algebra $A$ we must substitute axiom (3) by an analogous condition that does not mention the unit (but see Section \ref{section4.local.units} for the case of local units). In \cite{azevedo2021partial} the authors define partial coactions of multiplier Hopf algebras on regular algebras, and their solution to this problem is to replace the vector $\rho(1_A) \otimes 1_H$ by an idempotent multiplier. 
We also need to require some minimal regularity conditions from the nonunital algebra $A$. 

Recall that the right annihilator of an algebra $A$ is the ideal $r(A)=\{a\in A\,|\, ba=0,\,\forall b\in A\}$ and the left annihilator of $A$ is the ideal $l(A)=\{b\in A\,|\, ba=0,\,\forall a\in A\}$. If $M$ is a left $A$-module, its torsion submodule is $t_A(M)=\{m\in M\,|\,am=0, \, \forall a\in A\}$, and we say that $M$ is torsion free if $t_A(M)=0$. We can define torsion free right modules analogously.

We note that for a free left $A$-module we have that $t_A(A^{(\Lambda)} )= r(A)^{(\Lambda)}$.

\bigskip 

Now let $A$ be an associative algebra and consider a linear map $\rho:A\to A\otimes H$; we will use the Sweedler notation $\rho(a)=\sum a^{(0)}\otimes a^{(1)}$. Consider the linear map $\mu : (A \otimes H) \otimes A \to A \otimes H$ given by
\[
\mu\bigl((b\otimes h)\otimes a\bigr)=\sum ba^{(0)}\otimes ha^{(1)}.
\]
By analogy with the theory of partial $H$-coactions on unital algebras, in what follows we derive necessary and sufficient conditions for $\mu$ to define a right $A$-module structure on $A \otimes H$, and for $I\otimes \varepsilon:A\otimes H\to A$ to be a right $A$-linear map. The next two lemmas correspond to Lemma~2.1 and item~(3) of Lemma~2.2 of \cite{caenepeel}.

\begin{lemma}\label{lema1}
Let $A$ be an associative algebra, let $H$ be a Hopf algebra and assume that $r(A)=0$. Consider a linear map $\rho:A\to A\otimes H, \ \rho(a)=\sum a^{(0)}\otimes a^{(1)}$.  The tensor product $A\otimes H$ is an $A$-bimodule via $a(b\otimes h)c=\sum abc^{(0)}\otimes hc^{(1)}$ if and only if $\rho(ab)=\rho(a)\rho(b)$, i.e., $\sum (ab)^{(0)}\otimes (ab)^{(1)}=\sum a^{(0)}b^{(0)}\otimes a^{(1)}b^{(1)}$.
\end{lemma}
\begin{proof}
    It is clear that we only have to check the right $A$-module structure. 
    Since $r(A) = 0$ and the left $A$-module $A \otimes H$ is a free $A$-module, it follows that  ${t_A(A\otimes H)=0}$. 
    Therefore
    \begin{eqnarray*} ((a\otimes h)b)c&=&(a\otimes h)(bc)\\
    &\Updownarrow&\\
 \sum ab^{(0)}c^{(0)}\otimes hb^{(1)}c^{(1)}&=&\sum a(bc)^{(0)}\otimes h(bc)^{(1)}.    
    \end{eqnarray*}
Then, if $h=1$ and $t_A(A\otimes H)=0$, we have that $\sum b^{(0)}c^{(0)}\otimes b^{(1)}c^{(1)}=\sum (bc)^{(0)}\otimes (bc)^{(1)}$, i.e., $\rho(b)\rho(c)=\rho(bc)$.
\end{proof}

\begin{lemma}\label{lema2}
Let $H$ be a Hopf algebra, let $A$ be an associative algebra, and assume that $r(A)=0$. The map $I\otimes \varepsilon:A\otimes H\to A$ is right $A$-linear if and only if $(I\otimes \varepsilon)\rho(a)=a$.
\end{lemma}
\begin{proof}
    Note that
    \begin{eqnarray*}
        (I\otimes \varepsilon)(a\otimes h)b&=&(I\otimes\varepsilon)(\sum ab^{(0)}\otimes hb^{(1)})\\
        &\Updownarrow&\\
        ab\varepsilon(h)&=&\sum ab^{(0)}\varepsilon(hb^{(1)})
    \end{eqnarray*}
    Then, if $h=1$ and $r(A)=0$, we have that $b=\sum b^{(0)}\varepsilon(b^{(1)})=(I\otimes \varepsilon)\rho(b)$.
\end{proof}

\begin{definition} \label{def:partial.coaction}Let $H$ be a Hopf algebra and let $A$ be an associative algebra.
A linear map ${\rho: A\longrightarrow A \otimes H}$ is a \emph{partial $H$-coaction} on $A$ if it satisfies, for every $a,b \in A$,
\begin{enumerate}[\normalfont(i)]
	\item $(I\otimes\varepsilon_H)\rho(a)=a$;
	\item $(\rho\otimes I)((b\otimes 1_H)\rho(a))=(\rho(b)\otimes 1_H)(I\otimes \Delta_H)\rho(a)$,
\end{enumerate}
In this case, $A$ will be called a \emph{partial $H$-comodule algebra}.

If, additionally, $\rho$ also satisfies the equality
\begin{enumerate}[\normalfont(i)]
\addtocounter{enumi}{2}
\item
$(\rho\otimes I)(\rho(a)(b\otimes 1_H))=[(I\otimes \Delta_H)\rho(a)](\rho(b)\otimes 1_H)$ 
\end{enumerate}
for every $a,b\in A$, then $\rho$ is a \emph{symmetrical} partial coaction.
\end{definition}
Note that we did not ask that $\rho$ be a multiplicative map. In fact, this is a consequence of the axioms of a partial coaction. 
\begin{proposition} \label{prop:multiplicative}
    Let $A = (A,\rho)$ be a partial $H$-comodule algebra.  Then $\rho(ab) = \rho(a) \rho(b)$ for all $a,b$ in $A$. 
\end{proposition}
\begin{proof}
If $A$ is a partial $H$-comodule algebra, then we have that for every $a,b\in A$,
\begin{eqnarray*}
\sum (ab^{(0)})^{(0)}\otimes (ab^{(0)})^{(1)}\otimes b^{(1)}&=& \sum a^{(0)}b^{(0)}\otimes a^{(1)}(b^{(1)})_{(1)}\otimes (b^{(1)})_{(2)}\\
&\Downarrow^{I\otimes I\otimes \varepsilon} &\\
\sum (ab^{(0)})^{(0)}(ab^{(0)})^{(1)}\varepsilon(b^{(1)})&=& \sum a^{(0)}b^{(0)}\otimes a^{(1)}(b^{(1)})_{(1)}\varepsilon((b^{(1)})_{(2)})\\
&\Downarrow &\\
\sum (ab)^{(0)}\otimes (ab)^{(1)}&=&\sum a^{(0)}b^{(0)}\otimes a^{(1)}b^{(1)},
\end{eqnarray*}
i.e., 
\begin{equation*}
    \rho(ab)=\rho(a)\rho(b).
\end{equation*}     
\end{proof}
In particular, if $A$ is a unital algebra endowed with a partial $H$-comodule algebra structure (according to Definition \ref{def:partial.coaction}), then $A$ is a partial $H$-comodule algebra in the usual sense (Definition \ref{def:partial.coaction.unital}).

As a consequence of Proposition \ref{prop:multiplicative}, we may rewrite axioms (ii) and (iii) of Definition \ref{def:partial.coaction} in Sweedler notation as
\begin{enumerate}
    \item[(ii')] $\sum b^{(0)}a^{(0)} \otimes b^{(1)}a^{(1)} \otimes a^{(2)} =
    \sum b^{(0)}a^{(0)} \otimes b^{(1)}(a^{(1)})_{(1)} \otimes (a^{(1)})_{(2)}$,
    \item[(iii')] $\sum a^{(0)} b^{(0)}\otimes a^{(1)}b^{(1)} \otimes a^{(2)} = \sum
    a^{(0)}b^{(0)} \otimes (a^{(1)})_{(1)} b^{(1)} \otimes (a^{(1)})_{(2)}$.
\end{enumerate}

\begin{remark}
Both \cite{timmermann}, where partial coactions of $C^*$-bialgebras were introduced, and \cite{azevedo2021partial}, where partial coactions of multiplier Hopf algebras were introduced, present versions of axiom (ii) in which an idempotent multiplier takes the place of the term $\rho(1_A) \otimes 1_H$ that appears in the definition of partial coactions on unital algebras (Definition \ref{def:partial.coaction.unital}).
\end{remark}

\begin{definition}[\cite{ferrazza1}]\label{partial.action}
Let $A$ be an associative algebra. A linear map $\cdot:H\otimes A \longrightarrow A$, $h\otimes a \longmapsto h \cdot a$, is called a \textit{partial action} of $H$ on $A$ if, for all $a,b \in A$ and $h,k \in H$,
\begin{enumerate}[\normalfont(i)]
    \item $1_H\cdot a=a$;
    \item $h\cdot (a(k\cdot b))=\sum(h_{(1)}\cdot a)(h_{(2)}k\cdot b)$.
\end{enumerate}
In this case, $A$ is called a \textit{partial $H$-module algebra}. We say that the partial action is \textit{symmetrical} if, additionally,
\begin{enumerate}[\normalfont(i)]
\addtocounter{enumi}{2}
    \item $h\cdot ((k\cdot b)a)=\sum(h_{(1)}k\cdot b)(h_{(2)}\cdot a)$.
\end{enumerate}
for all $a,b\in A$ and $h,k\in H$.
\end{definition}

It is well known that when $A$ is a unital algebra and $H$ is finite-dimensional, the linear map $\rightharpoonup :H^*\otimes A \longrightarrow A$ that sends $f\otimes a$ to $f\rightharpoonup a=\sum a^{(0)}f(a^{(1)})$ is a partial action if and only if $\rho$ is a partial coaction. The same holds true for associative algebras.

\begin{proposition} \label{prop.fin.dim}
Let $H$ be a finite-dimensional Hopf algebra.
\begin{enumerate}[\normalfont(1)]
\item If $A$ is a partial $H$-comodule algebra with (symmetrical) partial coaction $\rho$, then $A$ is a partial $H^*$-module algebra with (symmetrical) partial action $\rightharpoonup: H^*\otimes A\longrightarrow A$ given by $f\rightharpoonup a=\sum a^{(0)}f(a^{(1)})$;
\item If $A$ is a partial $H$-module algebra with (symmetrical) partial action $\cdot$, then $A$ is a partial $H^*$-comodule algebra with (symmetrical) partial coaction $\rho: A\longrightarrow A\otimes H^*$ given by $\rho(a)=\sum_i d_i\cdot a\otimes p_i$, where $\dim H=n$, $\{d_i\}_{i=1}^n$ is a basis for $H$ and $\{p_i\}_{i=1}^n$ is its dual basis in $H^*$.
\end{enumerate}
\end{proposition}
The proof of the above result follows directly from the definitions. 

In \cite{alves2} it was proved that when $H$ is a cocommutative Hopf algebra the tensor product of two partial $H$-module categories is also an $H$-module category. There are analogous results for partial coactions and actions on associative algebras, as we show next.

\begin{theorem} \label{thm:commutative.tensor}Let $H$ be a commutative Hopf algebra. If $A$ and $B$ are both (symmetrical) partial $H$-comodule algebras, then $A\otimes B$ is a (symmetrical) partial $H$-comodule algebra via $\rho(a\otimes b)=\sum a^{(0)}\otimes b^{(0)}\otimes a^{(1)}b^{(1)}$.
\end{theorem} 
\begin{proof}
We ask for $H$ to be commutative because, for every $a, b\in A$, $x, y\in B$,
\begin{eqnarray*}
    (\rho\otimes I)((a\otimes x\otimes 1_H)\!\!\!\!\!\!\!&&\!\!\!\!\!\!\!\rho(b\otimes y))\\&=&(\rho\otimes I)\sum (ab^{(0)}\otimes xy^{(0)}\otimes b^{(1)}y^{(1)})\\
    &=&\sum (ab^{(0)})^{(0)}\otimes (xy^{(0)})^{(0)}\otimes (ab^{(0)})^{(1)}(xy^{(0)})^{(1)}\otimes b^{(1)}y^{(1)}\\
    &=& \sum a^{(0)}b^{(0)}\otimes x^{(0)}y^{(0)}\otimes a^{(1)}\underbrace{(b^{(1)})_{(1)}x^{(1)}}(y^{(1)})_{(1)}\otimes (b^{(1)})_{(2)}(y^{(1)})_{(2)}
\end{eqnarray*}
and
\begin{eqnarray*}
    (\rho(a\otimes x)\otimes 1_H)\!\!\!\!\!\!\!&&\!\!\!\!\!\!\!(I\otimes \Delta)\rho(b\otimes y)\\&=&\sum(a^{(0)}\otimes x^{(0)}\otimes a^{(1)}x^{(1)}\otimes 1_H)(b^{(0)} \otimes y^{(0)}\otimes (b^{(1)}y^{(1)})_{(1)}\otimes (b^{(1)}y^{(1)})_{(2)})\\
    &=&\sum a^{(0)}b^{(0)} \otimes x^{(0)}y^{(0)}\otimes a^{(1)}x^{(1)}(b^{(1)}y^{(1)})_{(1)}\otimes (b^{(1)}y^{(1)})_{(2)}\\
    &=&\sum a^{(0)}b^{(0)} \otimes x^{(0)}y^{(0)}\otimes a^{(1)}\underbrace{x^{(1)}(b^{(1)})_{(1)}}(y^{(1)})_{(1)}\otimes (b^{(1)})_{(2)}(y^{(1)})_{(2)}
\end{eqnarray*}
hence if $H$ is commutative, then $(\rho\otimes I)((a\otimes x\otimes 1_H)\rho(b\otimes y))=(\rho(a\otimes x)\otimes 1_H)(I\otimes \Delta)\rho(b\otimes y)$, i.e., item (ii) of Definition \ref{def:partial.coaction} holds true.
\end{proof}

\begin{corollary}Let $H$ be a commutative Hopf algebra. If $A$ and $B$ are (symmetrical) partial $H$-comodule unital algebras, then $A\otimes B$ is a (symmetrical) partial $H$-comodule unital algebra.
\end{corollary}

\begin{theorem} 
    Let $H$ be a cocommutative Hopf algebra. If $A$ and $B$ are both (symmetrical) partial $H$-module algebras (not necessarily unital), then $A\otimes B$ is a (symmetrical) partial $H$-module algebra via 
    \[h\cdot (a \otimes b)=\sum h_{(1)}\cdot a\otimes h_{(2)}\cdot b.
    \]
\end{theorem}
The proof of this result is analogous to the one of Theorem \ref{thm:commutative.tensor} and thus it will be omitted. 

%%%%%%%% Globalização de uma coação %%%%%%%%%%%%%

\section{Globalization of a partial coaction\label{section3:globalization}}

A source of partial coactions on unital algebras is the construction of the induced coaction. If $B$ is an $H$-comodule algebra with coaction $\gamma: B \to B \otimes H$ and $e \in B$ is a central idempotent, then the unital ideal $A = eB$ is a partial $H$-comodule (unital) algebra with coaction $\rho(b) = (e \otimes 1_H) \gamma (b)$.  In \cite{alves3}, Alves and Batista proved the globalization theorem which says that every partial coaction on a unital algebra is obtained in this manner, that is, they proved that there exists an \emph{enveloping coaction} for each partial coaction.   

\begin{definition}[\cite{alves3}] \label{def:globalization.unital.coaction}
Let $H$ be a Hopf algebra and $A$ a partial $H$-comodule unital algebra with coaction $\rho$. A pair $(B,\theta)$ is called an enveloping coaction (or globalization) for $\rho$ if:
\begin{enumerate}[\normalfont(i)]
	\item $B$ is an $H$-comodule algebra with coaction $\gamma$;
	\item $\theta: A \longrightarrow B$ is a monomorphism of algebras;
	\item $\theta(A)$ is an ideal of $B$;
	\item $(\theta \otimes I)\rho(a)=(\theta(1_A)\otimes 1_H)\gamma(\theta(a))$, for all $a\in A$;
	\item $B$ is generated by $\theta(A)$ as an $H$-comodule.
\end{enumerate}  
\end{definition}
 Items (iv) and (v) imply that $\theta(1_A)$ is a central idempotent of $B$, and then item (iv) has the following interpretation: $\theta$ is a morphism of partial coactions, where $\theta(A)$ has the induced partial coaction $\overline{\gamma}(\theta(a))=(\theta(1_A)\otimes 1_H)\gamma(\theta(a))$. 
%%%%%%%%%%%%%%%%%%%%%%%%%%%%%%%%%%%%%%%%

In this section we show that every partial $H$-coaction on an associative algebra has a globalization (Theorem \ref{theorem.globalization}). However, not only does item (iv) of Definition \ref{def:globalization.unital.coaction} have to be modified, but the absence of a unit also gets in the way of viewing the original partial coaction as an induced one. Nevertheless, for algebras with local units there is a notion of induced partial coaction (Proposition \ref{coacparcialideal}), and every partial coaction is induced (Proposition \ref{prop:glob.local.units}).

We also introduce partial $H$-comodule $\K$-linear categories. We prove a globalization theorem for partial coactions on categories, and we relate partial $H$-comodule categories with a class of partial $H$-comodule algebras with local units.

\subsection{Globalization of partial coactions on associative algebras}

In \cite{ferrazza1}, the authors investigated partial actions of Hopf algebras on algebras without identity and defined a globalization of a partial action. Also, they proved that when there exists a globalization $(B,\theta)$ of a partial action $\cdot:H\otimes A\to A$ and the right (or left) annihilator of the nonunital algebra $A$ was trivial, then the partial action on $\theta(A)$ was induced by the action $\triangleright:H\otimes B\to B$ through the projection $\pi(h\triangleright\theta(a))=\theta(h\cdot a)$. 
We will discuss briefly the definition of globalization for a partial action.

\begin{definition}[\cite{ferrazza1}] A \textit{globalization} for the partial action $\cdot$ is a pair $(B,\theta)$ such that 
	\begin{enumerate}[\normalfont(i)]
		\item $B$ is an $H$-module algebra (possibly nonunital), with action $\triangleright$;
		\item $\theta: A \longrightarrow B$ is a monomorphism of algebras;
		\item For every $a,b\in A$ and $h\in H$, we have $\theta((h\cdot a)b)=(h\triangleright \theta(a))\theta(b)$ and $\theta(b(h\cdot a))=\theta(b)(h\triangleright \theta(a))$; 
		\item $B=H\triangleright\theta(A)$.
	\end{enumerate}
\end{definition}
In \cite{ferrazza1} we also asked that 
\begin{enumerate}[\normalfont(i)]
\addtocounter{enumi}{4}
    \item 	 $h\bullet\theta(a)=\theta(h\cdot a)$ defines a partial action of $H$ on $\theta(A)$.
\end{enumerate}

However, (v) follows from items (i)-(iv): 
the map $\bullet:H\otimes \theta(A)\to \theta(A)$ given by $h\bullet\theta(a)=\theta(h\cdot a)$ is always well defined and is a partial action. The fact that it is a partial action part is a consequence of items (i) and (ii), and the fact that it is well defined is assured by the equality $\bullet=\theta\circ\cdot\circ (I\otimes \theta^{-1})$, which can be illustrated by the following commutative diagram.

\begin{center}
	$\xymatrixcolsep{3pc}\xymatrix{
		H\otimes \theta(A) \ar[r]^{\bullet}\ar[d]_{I\otimes \theta^{-1}} & \theta(A)\\
		H\otimes A\ar[r]_{\cdot} & A\ar[u]^{\theta}}$
\end{center}

On the other hand, we note that items (iii) and (iv) already imply that $\theta(A)$ is an ideal of $B$.
Dualizing items (i)-(iv), and adding the condition that $\theta(A)$ be an ideal of $B$ we obtain the definition for the globalization of a  partial $H$-coaction on a nonunital algebra $A$. 

\begin{definition}\label{def.globalization} Let $A$ be a partial $H$-comodule algebra with symmetrical partial coaction $\rho$. A pair $(B,\theta)$ will be called a globalization for the partial coaction $\rho$ if
\begin{enumerate}[\normalfont(i)]
	\item $B$ is an $H$-comodule algebra with coaction $\gamma$;
	\item $\theta: A\longrightarrow B$ is a monomorphism of algebras;
	\item $\theta(A)$ is an ideal of $B$;
	\item $(\theta\otimes I)((b\otimes 1_H)\rho(a))=(\theta(b)\otimes 1_H)\gamma(\theta(a))$, \\ $(\theta\otimes I)(\rho(a)(b\otimes 1_H))=\gamma(\theta(a))(\theta(b)\otimes 1_H)$, for every $a,b\in A$;
	\item $B$ is generated by $\theta(A)$ as an $H$-comodule.
\end{enumerate}
\end{definition}
Note that $(B,\theta)$ is a globalization in the sense of \cite{alves3} when $A$ is a unital algebra. Also, observe that if $(B,\theta)$ is a globalization of the partial coaction $\rho:A\to A\otimes H$, $\rho(a)=\sum a^{(0)}\otimes a^{(1)}$, then the map $\Phi:\theta(A)\to \theta(A)\otimes H$ given by $\Phi(\theta(a))=\sum \theta(a^{(0)})\otimes a^{(1)}=(\theta\otimes I)\rho(a)$ defines a partial coaction. The map $\Phi$ is well defined because $\Phi=(\theta\otimes I)\circ \rho\circ \theta^{-1}$ and $\theta$ is a monomorphism of algebras; moreover, it makes $\theta(A)$ a partial $H$-comodule algebra as a consequence of $\rho$ being a partial coaction.

We can illustrate how the partial coaction $\Phi$ works by the following commutative diagram, which shows, also, that $\theta$ is a morphism of partial coactions.

\begin{center}
	$\xymatrixcolsep{3pc}\xymatrix{
		A \ar[d]_{\rho} & \theta(A)\ar[l]_{\theta^{-1}}\ar[d]^{\Phi}\\
		A\otimes H\ar[r]_{\theta \otimes I} & \theta(A)\otimes H}$
\end{center}

\begin{remark}
If $(B,\theta)$ is a globalization of a partial coaction $\rho:A\to A\otimes H$, where $B$ is an $H$-comodule algebra with coaction $\gamma$ and $\Phi(\theta(a))=\sum \theta(a^{(0)})\otimes a^{(1)}$, then, from item (iv) of Definition \ref{def.globalization}, we have that
\[
(\theta(b)\otimes 1_H)\Phi(\theta(a)) = (\theta(b)\otimes 1_H)\gamma(\theta(a))
\]
since
\[
(\theta(b)\otimes 1_H)\Phi(\theta(a)) = 
(\theta(b)\otimes 1_H)((\theta\otimes I)\rho(a)) = 
(\theta\otimes I)((b\otimes 1_H)\rho(a))=(\theta(b)\otimes 1_H)\gamma (\theta(a)),
\]
% \begin{eqnarray*}
%     (\theta\otimes I)((b\otimes 1_H)\rho(a))&=&(\theta(b)\otimes 1_H)\gamma (\theta(a))\\
%     &\Downarrow&\\
%     (\theta(b)\otimes 1_H)((\theta\otimes I)\rho(a))&=&(\theta(b)\otimes 1_H)\gamma(\theta(a))\\
%     &\Downarrow&\\
%     (\theta(b)\otimes 1_H)\Phi(\theta(a))&=&(\theta(b)\otimes 1_H)\gamma(\theta(a))
% \end{eqnarray*}
which suggests that $\Phi$ and $\gamma$ are related in some way. In the unital case, taking $b=1_A$ and using $(1_A\otimes 1_H)\rho(a)=\rho(a)$, we have
\[
\Phi(\theta(a))=(\theta(1_A)\otimes 1_H)\gamma(\theta(a)).
\]
We will come back to this point in the case of algebras with local units (see Propositions \ref{coacparcialideal} and \ref{prop:glob.local.units}), where $\Phi$ becomes an induced partial coaction on $\theta(A)$.
\end{remark}

As in the unital case, the algebra $A\otimes H$ has the trivial $H$-comodule algebra structure given by $\delta=I\otimes \Delta$. Let $B$ be the subcomodule algebra of $A\otimes H$ generated by $\rho(A)$. Alternatively, $B$ is generated as an algebra by the vector subspace $H^*\rightharpoonup \rho(A)$, where $h^*\rightharpoonup \rho(a)=(I\otimes h^*)\delta(\rho(a))$ and $$\delta(\rho(a))=(I\otimes \Delta)(\sum a^{(0)}\otimes a^{(1)})=\sum a^{(0)}\otimes (a^{(1)})_{(1)}\otimes (a^{(1)})_{(2)} \ ;$$ in other words, $h^*\rightharpoonup \rho(a)=\sum a^{(0)}\otimes (a^{(1)})_{(1)}h^*((a^{(1)})_{(2)})$.

\begin{theorem}\label{theorem.globalization}
Let $\rho: A\to A\otimes H$ be a symmetrical partial coaction and let $B$ be the subcomodule algebra of $A\otimes H$ generated by $\rho(A)$. Then the pair $(B,\rho)$ is a globalization for $\rho$.
\end{theorem}
\begin{proof} Items (i) and (v) hold true by construction. Item (ii) holds since a partial coaction $\rho: A \to A \otimes H$ is a monomorphism of algebras, since $a=(I\otimes \varepsilon_H)\rho(a)$ and $\rho(ab)=\rho(a)\rho(b)$ for every $a,b \in A$.
For item (iii), if $a,b\in A$ and $h^*\in H^*$ then
\begin{eqnarray*}
\rho(b)(h^*\rightharpoonup \rho(a)) & = & 
\rho(b)(\sum a^{(0)}\otimes (a^{(1)})_{(1)}h^*((a^{(1)})_{(2)}))\\
&=&\sum b^{(0)}a^{(0)}\otimes b^{(1)}(a^{(1)})_{(1)}h^*((a^{(1)})_{(2)})\\
& = & (I\otimes I\otimes h^*) (\rho(b) \otimes 1_H )(I\otimes\Delta)\rho(a)\\
& = &  (I\otimes I\otimes h^*)(\rho\otimes I)((b\otimes 1_H)\rho(a))\\
&=&\sum b^{(0)}a^{(0)(0)}\otimes b^{(1)}a^{(0)(1)}h^*(a^{(1)})\\ 
& = & \rho(\sum ba^{(0)}h^*(a^{(1)})) \in \rho(A). 
\end{eqnarray*}
Hence $\rho(A)$ is a right ideal of $B$. To show that $\rho(A)$ is also a left ideal, we use an analogous computation and the hypothesis that $\rho$ is symmetrical.
Finally, item (iv) is an immediate consequence of the fact that $\rho$ is a symmetrical partial coaction.
\end{proof}

\begin{remark}[see \cite{alves3}]
When $H$ is a finite-dimensional Hopf algebra, there is another way to construct a globalization of a partial coaction $\rho:A\to A\otimes H$: consider the corresponding partial $H^*$-action and its minimal globalization $B=H^*\triangleright \varphi(A)$, where $\varphi: A \to \Hom(H^*, A)$ is given by $\varphi(a)(h^*)=h^*\cdot a=a^{(0)}h^*(a^{(1)})$ and
\[(k^*\triangleright \varphi(a))(h^*)=h^*k^*\cdot a=a^{(0)}(h^*k^*)(a^{(1)})=a^{(0)}h^*((a^{(1)})_{(1)})k^*((a^{(1)})_{(2)}).\]
Hence $B$ is an $H$-comodule algebra with $\gamma:B\to B\otimes H$, given by $\gamma(k^*\triangleright\varphi(a))=\sum_i h_i^*k^*\triangleright \varphi(a)\otimes h_i$, for a chosen basis $\{h_i\}$ of $H$. In this case, $(B,\varphi)$ is also a globalization of $\rho$.
\end{remark}

\subsection{The minimal globalization of a partial action by a finite dimensional Hopf algebra}

In \cite{alves3} there is a description of the globalization of a partial coaction via the globalization of the induced partial action. Here, we will prove that when $H$ is finite dimensional, the minimal globalization of a partial action coincides with a globalization of the induced partial coaction.

\begin{proposition}\label{prop2}
Let $H$ be a finite-dimensional Hopf algebra, $A$ an associative algebra, and $\cdot: H\otimes A\to A$ a symmetrical partial action. Let $\rho: A\to A\otimes H^*$ be the induced partial coaction $\rho(a)=\sum_{i=1}^{n}h_i\cdot a\otimes h_i^*$, where $\{h_i\}_{i=1}^{n}$ is a basis of $H$ and $\{h_i^*\}_{i=1}^{n}$ is the dual basis, and consider the action $\rightharpoonup: H\otimes A\otimes H^*\to A\otimes H^*$ given by $h\rightharpoonup (a\otimes k^*)=\sum k^*_{(2)}(h)a\otimes k^*_{(1)}$. Then $(H\rightharpoonup \rho(A),\rho)$ is the minimal globalization of $\cdot: H\otimes A\to A$.
\end{proposition}
\begin{proof}First, let us prove that $B=H\rightharpoonup\rho(A)$ is an $H$-module algebra. For this, we only need to prove that
\begin{eqnarray*}
(h_k\rightharpoonup \rho(a))(h_l\rightharpoonup \rho(b))=[\sum_{i=1}^{n}(h^*_i)_{(2)}(h_k)h_i\cdot a\otimes (h^*_i)_{(1)}][\sum_{j=1}^{n}(h^*_j)_{(2)}(h_l)h_j\cdot b\otimes (h^*_j)_{(1)}] \in B,
\end{eqnarray*} 
for every $k,l\in\{1,2,\cdots,n\}$, $a, b\in A$. In fact, for any $h\in H$,
\begin{eqnarray*}
(I\otimes h)\!\!\!\!\!\!\!\!\!\!&&\!\!\!\!\![\sum_{i=1}^{n}(h^*_i)_{(2)}(h_k)h_i\cdot a\otimes (h^*_i)_{(1)}][\sum_{j=1}^{n}(h^*_j)_{(2)}(h_l)h_j\cdot b\otimes (h^*_j)_{(1)}]\\
&=&(I\otimes h)[\sum_{i,j}(h^*_i)_{(2)}(h_k)(h^*_j)_{(2)}(h_l)(h_i\cdot a)(h_j\cdot b)\otimes (h^*_i)_{(1)}(h^*_j)_{(1)}]\\
&=&\sum_{i,j,(h)}(h^*_i)_{(2)}(h_k)(h^*_j)_{(2)}(h_l)(h_i\cdot a)(h_j\cdot b)\otimes (h^*_i)_{(1)}(h_{(1)})(h^*_j)_{(1)}(h_{(2)})\\
&=&\sum_{i,j,(h)}(h^*_i)(h_{(1)}h_k)(h^*_j)(h_{(2)}h_l)(h_i\cdot a)(h_j\cdot b)\\
&=&\sum (h_{(1)}h_k\cdot a)(h_{(2)}h_l\cdot b)\\
&=&\sum (h_{(1)}(h_k)_{(1)}\cdot a)(h_{(2)}(h_k)_{(2)}S((h_k)_{(3)})h_l\cdot b)\\
&=&\sum(h(h_k)_{(1)}\cdot a(S((h_k)_{(2)})h_l\cdot b))\\
&=&\sum_i h^*_i(h(h_k)_{(1)})h_i\cdot a(S((h_k)_{(2)})h_l\cdot b))\\
&=&\sum_i (h^*_i)_{(1)}(h)(h^*_i)_{(2)}((h_k)_{(1)})h_i\cdot a(S((h_k)_{(2)})h_l\cdot b))\\
&=&(I\otimes h)(\sum_i (h^*_i)_{(2)}((h_k)_{(1)})h_i\cdot a(S((h_k)_{(2)})h_l\cdot b))\otimes (h^*_i)_{(1)})\\
&=&(I\otimes h)(\sum (h_k)_{(1)}\rightharpoonup \rho (a(S((h_k)_{(2)})h_l\cdot b))).
\end{eqnarray*}
Hence
\[
\Bigl[\sum_{i=1}^{n}(h^*_i)_{(2)}(h_k)h_i\cdot a\otimes (h^*_i)_{(1)}\Bigr]
\Bigl[\sum_{j=1}^{n}(h^*_j)_{(2)}(h_l)h_j\cdot b\otimes (h^*_j)_{(1)}\Bigr]
=\sum (h_k)_{(1)}\rightharpoonup \rho\bigl(a(S((h_k)_{(2)})h_l\cdot b)\bigr)\in B.
\]
The other requisites, except for minimality, follow from
\begin{eqnarray*}
\rho(a(k\cdot b))&=&\rho(a)(k\rightharpoonup\rho(b))\\
\rho((k\cdot b)a)&=&(k\rightharpoonup\rho(b))\rho(a).
\end{eqnarray*}
Finally, for the minimality, suppose that $\sum_i kh^i\cdot a^i=0$ for all $k\in H$. Then, for every $k\in H$
\begin{eqnarray*}
(I \otimes k)(\sum_i h^i\rightharpoonup \rho(a^i))&=&(I\otimes k)(\sum_{i,j}h^i\rightharpoonup h_j\cdot a^i\otimes h^*_j)\\
&=&(I\otimes k)(\sum_{i,j} (h^*_j)_{(2)}(h^i)h_j\cdot a^i\otimes (h^*_j)_{(1)})\\
&=&\sum_{i,j}h^*_j(kh^i)h_j\cdot a^i\\
&=&\sum_i kh^i\cdot a^i=0,
\end{eqnarray*}
hence $\sum_i h^i\rightharpoonup \rho(a^i)=0$.
\end{proof} 

This proposition can be illustrated by the following commutative diagram:
\[
\xymatrix{H\otimes A\otimes H^* \ar[rr]^{I_H\otimes \theta} \ar[d]_{\rightharpoonup} & {} & H\otimes \Hom(H,A) \ar[d]^{\triangleright} \\
A\otimes H^* \ar[rr]^{\theta} & {} & \Hom(H,A) \\
{} & A \ar[lu]^{\rho} \ar[ru]_{\varphi} & {}}
\]
where $\theta$ is the natural isomorphism.

\begin{definition}[\cite{ferrazza1}]
Let $G$ be a finite group. A good partial $G$-grading of $FMat_{\mathbb{N}}(\Bbbk)$ is a partial action of $(\Bbbk G)^*$ on $FMat_{\mathbb{N}}(\Bbbk)$ such that each elementary matrix $E_{ij}$ is a simultaneous eigenvector for all elements of the canonical basis $\{p_g\,;\, g \in G\}$ of $(\Bbbk G)^*$.
\end{definition}

\begin{example}
In \cite{ferrazza1}, it was proved that if $G$ is a finite abelian group such that $\operatorname{char}(\Bbbk)\nmid |G|$, $H=(\Bbbk G)^*$, $A=FMat_{\mathbb{N}}(\Bbbk)$, and $\cdot: H\otimes A\to A$ is a good partial $G$-grading on $A$, then $(B,\theta)$, where $B=\bigoplus_{i,j\in \mathbb{N}}B_{ij}$ is an $H$-module subalgebra of $FMat_{\mathbb{N}}(\Bbbk G)$, where $B_{ij}=\bigoplus_{g\in t_{ij}L} \Bbbk g$ and the $H$-action is induced by $p_g\rightharpoondown h=\delta_{g,h}h$, and $\theta : A \to B$ is given by
\[
E_{ij}\longmapsto \left(\delta_{(x,y),(i,j)}\frac{1}{|L|}\sum_{g\in t_{ij}L}g\right)_{x,y\in \mathbb{N}}=\left(\frac{1}{|L|}\sum_{g\in t_{ij}L}g\right)E_{ij},
\]
is a minimal globalization for the good partial $G$-grading on $A$.

Now, applying the previous proposition on the good partial $G$-grading on $A$, we have that the corresponding partial $\K G$-coaction is
\begin{eqnarray*}
\rho(E_{ij})&=&\sum_{g\in G}p_g\cdot E_{ij}\otimes g\\
&=&\sum_{g\in G}\frac{1}{|L|}\delta_{gt_{ij}L,t_{ij}L}E_{ij}\otimes g\\
&=&\sum_{g\in t_{ij}L}\frac{1}{|L|}E_{ij}\otimes g, 
\end{eqnarray*}
and by the algebra isomorphism $FMat_{\mathbb{N}}(\Bbbk )\otimes \Bbbk G\simeq FMat_\mathbb{N}(\Bbbk G)$, we may write 
\begin{eqnarray*}
\rho(E_{ij})=(\sum_{g\in t_{ij}L}\frac{1}{|L|}g)E_{ij}\in FMat_{\mathbb{N}}(\Bbbk G).
\end{eqnarray*}
This gives us the same minimal globalization presented in \cite{ferrazza1}.
\end{example}

\begin{example} Consider $A=\Bbbk $ and $G$ a finite group. The minimal globalization of a partial action $\cdot: \Bbbk G\otimes \Bbbk \to \Bbbk $ given by $g\cdot 1=\delta_{gH,H}$ where $H$ is a subgroup of $G$, is the pair $((\K H)^*,\theta)$, where $\theta(1)=\sum_{h\in H}p_h$.
\end{example}
\begin{example} Let $G$ be a finite group and $A$ a partially $G$-graded associative algebra, i.e., $A$ is a partial $(\K G)^*$-module algebra with symmetrical partial action. Note that for every $a\in A$, we have
\begin{eqnarray*}
\rho(a)=\sum_{g\in G}p_g\cdot a\otimes g,
\end{eqnarray*}
and
\begin{eqnarray*}
p_h\rightharpoonup\rho(a)=p_h\cdot a\otimes h.
\end{eqnarray*}
Then the minimal globalization of the original partial action is given by an algebra that is generated as vector space by the elements $p_g\cdot a\otimes g$, with $g\in G$ and $a\in A$. Note also that
\begin{eqnarray*}
p_{gh}\cdot(a(p_h\cdot b))&=&\sum_{t\in G}(p_{ght^{-1}}\cdot a)(p_tp_h\cdot b)\\
&=& (p_g\cdot a)(p_h\cdot b),
\end{eqnarray*}
then we can identify the minimal globalization as the $G$-graded algebra $B=\oplus_{g\in G}B_g$, where $B_g=p_g\cdot A$, and the monomorphism of algebras $\theta: A\to B$ is given by $\theta(a)=\sum_{g\in G}p_g\cdot a$.
\end{example}

Note that if we consider a partial $G$-grading on the algebra with local units $A=FMat_\mathbb{N}(\Bbbk )$, then its minimal globalization is $B=\oplus_{g\in G}(p_g\cdot FMat_\mathbb{N}(\Bbbk ))$, which can be identified as a subalgebra of $FMat_\mathbb{N}(\Bbbk G)$, even if the partial grading is not a good grading. This proves the last result of the present section.

\begin{proposition} The minimal globalization of any partial $G$-grading of $FMat_\mathbb{N}(\Bbbk )$ is a subalgebra of $FMat_\mathbb{N}(\Bbbk G)$.
\end{proposition}

\section{Partial $H$-comodule categories and partial $H$-comodule algebras with local units \label{section4.local.units}}

\subsection{Partial coactions on algebras with local units}
 
\begin{definition}
Let $A$ be a $\Bbbk$-algebra. 
	A \emph{system of local units} (s.l.u.) of $A$ is a subset $S \subset A$ such that $e^2=e$ for every $e \in S$ and, for every \emph{finite} subset $F$ of $A$, there exists $e \in S$ such that $e a=ae =a$ for every $a\in F$. In this case we say that $A$ is an \emph{algebra with local units}.
\end{definition}

The next result presents a slight reformulation of the definition of partial coaction when $A$ is an algebra with local units which  is more appropriate to this context.

\begin{proposition}\label{prop.partial.coaction.local.units} Let $H$ be a Hopf algebra, $A$ an algebra with local units and $S$ an s.l.u. for $A$. The linear map $\rho: A\longrightarrow A \otimes H$ is a partial coaction if and only if, for every $a,b \in A$,
\begin{enumerate}[\normalfont(1)]
	\item $(I\otimes\varepsilon_H)\rho(a)=a$;
	\item $\rho(ab)=\rho(a)\rho(b)$;
	\item $(\rho\otimes I)\rho(a)=(\rho(e)\otimes 1_H)(I\otimes \Delta_H)\rho(a)$, for every $e\in S$ such that $(e\otimes 1_H)\rho(a)=\rho(a)$.
\end{enumerate}
The partial coaction $\rho$ is symmetrical if and only if, additionally,
\[
(\rho\otimes I)\rho(a)=[(I\otimes\Delta_H)\rho(a)](\rho(f)\otimes 1_H)
\]
for every $f\in S$ such that $\rho(a)(f\otimes 1_H)=\rho(a)$.
\end{proposition}
\begin{proof}
    Assume that (1), (2), (3) hold. Given $a, b  \in A$, we may choose a local unit $e$ such that 
    \[
   \sum  ea^{(0)} \otimes a^{(1)}=\sum a^{(0)} \otimes a^{(1)} 
    \]
    and also $b = be$, hence, $\sum b^{(0)} \otimes b^{(1)} = \sum b^{(0)}e^{(0)} \otimes b^{(1)}e^{(1)}$. In Sweedler notation, (3) is  
    \begin{align*}
        \sum a^{(0)} \otimes a^{(1)}  \otimes a^{(2)} = \sum e^{(0)}a^{(0)} \otimes e^{(1)}(a^{(1)})_{(1)}  \otimes (a^{(2)})_{(1)};  
    \end{align*}
    multiplying this equality by $\rho(b) \otimes I_H$, 
    \begin{align*}
        \sum b^{(0)}a^{(0)} \otimes b^{(1)}a^{(1)}  \otimes a^{(2)} & = \sum b^{(0)}e^{(0)}a^{(0)} \otimes b^{(1)}e^{(1)}(a^{(1)})_{(1)}  \otimes (a^{(2)})_{(1)}  
        \\ 
        & = \sum 
b^{(0)}a^{(0)} \otimes b^{(1)}(a^{(1)})_{(1)}  \otimes (a^{(2)})_{(1)},    
    \end{align*}
    and thus item (ii) of Definition \ref{def:partial.coaction} holds. Conversely, given (ii) of Definition \ref{def:partial.coaction}, substituting $b=e$ we obtain condition (3) of this proposition; (1) coincides with (i), and (2) is a consequence of (i) and (ii) (see Proposition \ref{prop:multiplicative}).
\end{proof}

Note that, as mentioned in \cite{ferrazza1} for partial actions, if $A$ is a partial $H$-comodule algebra with local units satisfying conditions (1)--(3) for some s.l.u. $S$, then the same conditions hold for every s.l.u. $T$ of $A$.

%%%%%%%%%%%%%%%%%5 Partial coaction on categories %%%%%%%%%%%%%%%%

\subsection{Partial comodule categories}

In \cite{ferrazza1} it was shown that partial $H$-actions on algebras with local units are connected to partial $H$-actions on $\K$-linear categories. We will show that the same holds for partial coactions.

$H$-coactions on categories were inspired by previous results on group gradings of $\K$-linear categories.
 By a \emph{$\K$-linear category} we mean a $V$-category over the monoidal category $V$ of $\K$-vector spaces, that is, a category $\mC$ where each set of morphisms ${}_y\mC_x$ is a $\K$-vector space and each composition mapping ${}_z\mC_y \times {}_y\mC_x \to {}_z\mC_x$ is $\K$-bilinear \cite{cibils-solotar-galois}. We will often say only ``linear category''.
\begin{definition}[\cite{cibils2006skew}]
Let $G$ be a group and let $\mathcal{C}$ be a linear category. A $G$-grading of $\mathcal{C}$ is a direct sum decomposition of morphism spaces
\[
 {}_y\mathcal{C}_x = \bigoplus_{g\in G}  \ {{}_y\mathcal{C}_x}^g
\]
such that
\[
 {_z\mathcal{C}_y}^g \circ  {{}_y\mathcal{C}_x}^h \subseteq  {{}_z\mathcal{C}_x}^{gh}.
\]
\end{definition}

\begin{example}
Consider $A=Mat_{n\times n}(\Bbbk)$ with a good $G$-grading (see \cite{ion}), where $G$ is a finite group. Then, the linear category $\mathcal{C}$ with $n$ objects and $_i\mathcal{C}_j=\Bbbk$ is an $H$-comodule category with coaction $\rho_{(i,j)}(\alpha)=\alpha\otimes g$ if and only if $E_{ij}\in A_g$.
\end{example}

\begin{example}[\cite{cibils2006skew,cibils2022resolving}]
Let $\mathcal{C}$ be a linear category, and consider an action of a group $G$ on $\mathcal{C}$ by $\K$-linear automorphisms. The \emph{skew category} associated to the $G$-action is the linear category $\mathcal{C}[G]$ with the same objects as $\mathcal{C}$, with morphism spaces given by
\[
{}_y\mathcal{C}[G]_x = \bigoplus_{r \in G}  {}_y\mathcal{C}_{rx}.
\]
A $G$-grading is introduced by declaring each ${}_y\mathcal{C}_{rx}$ to be the component of degree $r$ of ${}_y\mathcal{C}[G]_x$, that is, ${{}_y\mathcal{C}[G]_x}^r = {}_y\mathcal{C}_{rx}$. The composition of morphisms is defined as follows: if $f \in {}_y\mathcal{C}[G]_x$ has degree $r$ and $g \in {}_z\mathcal{C}[G]_y$ has degree $s$, then $g \circ_{\mathcal{C}[G]} f = g \circ_{\mathcal{C}}(sf)$, which lands in ${}_z\mathcal{C}[G]_{srx}$ and hence has degree $sr$. It follows that each $G$-action on $\mathcal{C}$ by $\K$-linear automorphisms yields a $G$-graded category.
\end{example}

$H$-comodule categories were introduced in \cite{herscovich2007hochschild} modeled after $G$-gradings of categories. 

\begin{definition}
[\cite{herscovich2007hochschild}] Let $H$ be a Hopf algebra. 
A linear category $\mathcal{C}$ is an \textit{$H$-comodule category} if there exists a family of linear maps $\rho=\{\rho_{(x,y)}:\,_y\mathcal{C}_x\to\,_y\mathcal{C}_x\otimes H\}_{x,y\in\mathcal{C}_0}$ such that, for every $_yf_x\in\,_y\mathcal{C}_x$, $_zg_y\in\,_z\mathcal{C}_y$,
\begin{enumerate}[\normalfont(i)]
	\item $\rho_{(x,z)}(gf)=\rho_{(y,z)}(g)\rho_{(x,y)}(f)$;
	\item $\rho_{(x,x)}(\,_x1_x)=\,_x1_x\otimes 1_H$;
	\item $(I\otimes \Delta)\rho_{(x,y)}=(\rho_{(x,y)}\otimes I)\rho_{(x,y)}$;
	\item $(I\otimes \varepsilon)\rho_{(x,y)}(_yf_x)=\,_yf_x$.
\end{enumerate}
The family $\rho$ is the $H$-coaction on $\mathcal{C}$.
\end{definition}

\begin{example}
    Just as the case of $G$-graded algebras, a $G$-grading of $\mathcal{C}$ is the same as a structure of $\Bbbk G$-comodule category on $\mathcal{C}$. Given the $G$-grading one defines a $\K G$-coaction on homogeneous components by $\rho_{(x,y)} ({}_{y}f_{x})= {}_{y}f_{x} \otimes r$. Conversely, to each $\K G$-coaction on  $\mathcal{C}$ we assign a $G$-grading by declaring ${{}_{y}\mathcal{C}_x}^r$ to be the vector subspace of ${}_{y}\mathcal{C}_x$ consisting of the morphisms such that $\rho_{(x,y)} ({}_{y}f_{x})= {}_{y}f_{x} \otimes r$.
\end{example}
%%%%%%%%%%%%%%%%%%%%%%%%%%%%%%

\begin{definition}\label{def.partial.coaction.category} Let $H$ be a Hopf algebra. A linear category $\mathcal{C}$ is a \textit{partial $H$-comodule category} if there exists a family of linear maps $\rho=\{\rho_{(x,y)}:\,_y\mathcal{C}_x\to\,_y\mathcal{C}_x\otimes H\}_{x,y\in\mathcal{C}_0}$ such that, for every $_yf_x\in\,_y\mathcal{C}_x$, $_zg_y\in\,_z\mathcal{C}_y$,
\begin{enumerate}[\normalfont(i)]
	\item $(I\otimes \varepsilon)\rho_{(x,y)}(_yf_x)=\,_yf_x$
	\item $\rho_{(x,z)}(gf)=\rho_{(y,z)}(g)\rho_{(x,y)}(f)$;
	\item $(\rho_{(x,y)}\otimes I)\rho_{(x,y)}(f)=(\rho_{(y,y)}(\,_y1_y)\otimes 1_H)(I\otimes \Delta)\rho_{(x,y)}(f)$.
\end{enumerate}
In this case, $\rho$ will be called a \emph{partial $H$-coaction on $\mathcal{C}$}. The partial coaction $\rho$ will be called \emph{symmetrical} if, additionally, 
\begin{enumerate}[\normalfont(i)]
\addtocounter{enumi}{3}
    \item $(\rho_{(x,y)}\otimes I)\rho_{(x,y)}(f)=(I\otimes \Delta)\rho_{(x,y)}(f)[\rho_{(x,x)}(\,_x1_x)\otimes 1_H]$.
\end{enumerate}
\end{definition}

It is clear that every $H$-comodule category is a partial $H$-comodule category.

One can associate a partial $H$-comodule algebra with local units to each partial $H$-comodule category. 

\begin{definition}[\cite{ferrazza1}] Let $\mathcal{C}$ be a linear category. We will denote by $a(\mathcal{C})$ the algebra consisting of elements of the form $(\,_yf_x)_{x,y\in\mathcal{C}_0}$, with finite $_yf_x\neq 0$, where the multiplication is given by matrix multiplication and the composition of $\mathcal{C}$.
\end{definition}

\begin{remark}
    Note that $(\,_yf_x)_{x,y\in\mathcal{C}_0}$ is the matrix notation of the elements of $a(\mathcal{C})$, which, as vector space, is the direct sum $\bigoplus_{x,y\in \mathcal{C}_0} \,_y\mathcal{C}_x $. The algebra $a(\mathcal{C})$ has a canonical system of local units $S_0$ whose elements are finite sums of distinct elements $e_{xx}$, which are given by $_x1_x$ at the position $(x,x)$ and zero at every other position. Hence, for every finite subset $P\subset \mathcal{C}_0$, $P\neq \varnothing$, we may associate the local unit $e_P=\sum_{x\in P}e_{xx}$.
\end{remark}

\begin{proposition}
    Let $\mathcal{C}$ be a linear category. If $\mathcal{C}$ is a (symmetrical) partial $H$-comodule category with partial coaction $\{\rho_{(x,y)}:\,_y\mathcal{C}_x\to \,_y\mathcal{C}_x\otimes H\}_{x,y\in\mathcal{C}_0}$, then $a(\mathcal{C})$ is a (symmetrical) partial $H$-comodule algebra with local units with partial coaction defined by $\rho((\,_yf_x)_{x,y\in\mathcal{C}_0})=(\rho_{(x,y)}(\,_yf_x))_{x,y\in\mathcal{C}_0}$.
\end{proposition}
\begin{proof}
        Since $\rho_{(x,y)}(\,_y\mathcal{C}_x)\subseteq\,_y\mathcal{C}_x\otimes H$ and for a local unit $e_P$ we have
        \[
        (e_P\otimes 1_H)\rho((\,_yf_x)_{x,y\in\mathcal{C}_0})=\rho((\,_yf_x)_{x,y\in\mathcal{C}_0})
        \]
        if and only if $y\in P\subset \mathcal{C}_0$ for all $y$ such that $_yf_x\neq 0$, the calculations that prove that $\rho:a(\mathcal{C})\to a(\mathcal{C})\otimes H$ defines a (symmetrical) partial coaction on $a(\mathcal{C})$ follow directly from Definition \ref{def.partial.coaction.category}.
\end{proof}

It is possible to prove a converse of this result for a class of partial $H$-coactions. The definition of this class which is motivated by our previous work on partial actions, and we recall the corresponding concept of \emph{categorizable} partial action.  

\begin{definition}[\cite{ferrazza1}] If a linear map $\cdot : H\otimes A \longrightarrow A$ is a (symmetrical) partial action and there exists some s.l.u. $S$ such that $H\cdot e \subseteq eA e$, for every $e \in S$, we will say that this partial action is an \emph{$S$-categorizable (symmetrical) partial action}.
\end{definition}

Inspired by the definition above about categorizable partial actions, which induces partial actions on a category obtained from the respective algebra with local units (see \cite{ferrazza1}), we introduce the definition of \emph{categorizable partial coactions}.

\begin{definition} If $A$ is a partial $H$-comodule algebra with local units with (symmetrical) partial coaction $\rho$, and for some s.l.u. $S$ we have that $\rho(e)\in eAe\otimes H$ for every $e \in S$, then $\rho(eAf)\subseteq eAf\otimes H$ for every $e,f \in S$. In this case, $\rho$ will be called a (symmetrical) \textit{$S$-categorizable partial coaction}.
\end{definition}

\begin{definition}[\cite{ferrazza1}]
Let $A$ be an algebra and $S$ be a system of local units for $A$. The $\K$-linear category $\mathcal{C}^{S}(A)$ has $S$ as its set of objects. For any two objects $e,f \in S$ we associate the vector space of morphisms  $\mathcal{C}^{S}(A)(e,f)=\,_f\mathcal{C}^{S}(A)_e=fAe$, and composition is just the restriction of multiplication of $A$.  
\end{definition}

In \cite{ferrazza1}, it was proved that if $A$ is a partial $H$-module algebra with local units and with a (symmetrical) $S$-categorizable partial action, then the linear category $\mathcal{C}^S(A)$ is a partial $H$-module category. Analogously, we have the following result for (symmetrical) partial coactions.

\begin{proposition}
Let $A$ be an algebra with an s.l.u. $S$. If $A$ is a (symmetrical) partial $H$-comodule algebra with an $S$-categorizable (symmetrical) partial coaction $\rho:A\to A\otimes H$, then $\mathcal{C}^S(A)$ is a (symmetrical) partial $H$-comodule category with partial coaction $\{\rho_{(e,f)}=\rho|_{fAe}\}_{e,f\in S}$.
\end{proposition}

\begin{proof}
Given an $S$-categorizable (symmetrical) partial coaction $\rho:A\to A\otimes H$, since $\rho_{(e,f)}$ is the restriction of $\rho$ to $fAe$ and $\rho$ is a (symmetrical) partial coaction, the proof that $\{\rho_{(e,f)}\}_{e,f\in S}$ defines a (symmetrical) partial coaction on $\mathcal{C}^S(A)$ follows straightforwardly from Proposition \ref{prop.partial.coaction.local.units}.
\end{proof}

We close this subsection with two results analogous to Proposition \ref{prop.fin.dim} and Theorem \ref{thm:commutative.tensor}, respectively. The proofs are similar and thus will be omitted. 

\begin{proposition} Let $H$ be a finite dimensional Hopf algebra.
\begin{enumerate}[\normalfont(1)]
\item If $\mathcal{C}$ is a partial $H$-module category, then $\mathcal{C}$ is a partial $H^*$-comodule category;
\item If $\mathcal{C}$ is a partial $H$-comodule category, then $\mathcal{C}$ is a partial $H^*$-module category.
\end{enumerate} 
\end{proposition}
\begin{theorem}
Let $H$ be a commutative Hopf algebra. If $\mathcal{A}$ and $\mathcal{B}$ are both (symmetrical) partial $H$-comodule categories, then the tensor product category $\mathcal{A}\otimes \mathcal{B}$ is a (symmetrical) partial $H$-comodule category via 
\[
\rho_{(x,y)} ({}_ya_{x} \otimes {}_yb_{x}) = \sum
{}_ya_{x}^{(0)} \otimes {}_yb_{x}^{(0)} \otimes {}_ya_{x}^{(1)} {}_yb_{x}^{(1)}. 
\]
\end{theorem}

%%%% section: globalization partial coaction on category 

\subsection{Globalization of a partial coaction on a $\K$-linear category}

In this section we prove that there always exists a globalization of a symmetrical partial coaction on a linear category. The globalization occurs in a linear \emph{semicategory}, which is essentially a linear category which may not have identities for some of its objects. This follows the pattern of globalizations of partial actions and coactions on unital algebras: the globalization is usually a nonunital algebra.  First, let us recall the definitions of semicategory and semifunctor.

\begin{definition}[\cite{alves2} (semicategories and semifunctors)]
\begin{enumerate}[\normalfont(i)]
    \item A linear semicategory $\mathcal{C}$ consists of a class $\mathcal{C}_0$ (its objects) and, for each $x, y \in \mathcal{C}_0$, a $\Bbbk$-linear space ${}_y\mathcal{C}_x$, the space of morphisms from $x$ to $y$, such that $\mathcal{C}$ satisfies all the axioms of a linear category except for the existence of an identity for every object.
    \item A linear semifunctor from the linear semicategory $\mathcal{C}$ to the linear semicategory $\mathcal{D}$ is a family of $\Bbbk$-linear maps ${}_yF_x : {}_y \mathcal{C}_x \to {}_{F(y)} \mathcal{C}_{F(x)}$, where $(x,y) \in \mathcal{C}_0 \times \mathcal{C}_0$, such that ${}_zF_x ({}_zg_y \circ {}_yf_x) = {}_zF_y ({}_zg_y) \circ {}_yF_x({}_yf_x)$ for every $x,y,z \in \mathcal{C}_0$, ${}_yf_x \in {}_y\mathcal{C}_x$, and ${}_zg_y \in {}_z\mathcal{C}_y$. As usual, we will simply write $F$ instead of ${}_yF_x$.
    \item Let $\mathcal{C}$ be a linear semicategory.
    \begin{enumerate}[\normalfont(a)]
        \item A $\mathcal{C}_0$-semicategory is a semicategory $\mathcal{D}$ whose class of objects is $\mathcal{C}_0$.
        \item A $\mathcal{C}_0$-subsemicategory $\mathcal{E}$ of a $\mathcal{C}_0$-semicategory $\mathcal{D}$ is a subsemicategory of $\mathcal{D}$ such that $\mathcal{E}_0 = \mathcal{D}_0 = \mathcal{C}_0$.
        \item A $\mathcal{C}_0$-semifunctor is a functor between $\mathcal{C}_0$-semicategories that is the identity on the objects.
    \end{enumerate}
\end{enumerate}
\end{definition}

\begin{definition}Let $H$ be a Hopf algebra. A linear semicategory $\mathcal{H}$ is an \textit{$H$-comodule semicategory} if there exists a family of linear maps $\rho=\{\rho_{(x,y)}:\,_y\mathcal{H}_x\to\,_y\mathcal{H}_x\otimes H\}_{x,y\in\mathcal{H}_0}$ such that, for every $_yf_x\in\,_y\mathcal{H}_x$, $_zg_y\in\,_z\mathcal{H}_y$,
\begin{enumerate}[\normalfont(i)]
	\item $\rho_{(x,z)}(gf)=\rho_{(y,z)}(g)\rho_{(x,y)}(f)$;
	\item $(I\otimes \Delta)\rho_{(x,y)}=(\rho_{(x,y)}\otimes I)\rho_{(x,y)}$;
	\item $(I\otimes \varepsilon)\rho_{(x,y)}(f)=f$.
\end{enumerate}
\end{definition}
The globalization of a symmetric partial coaction on a linear category $\mathcal{C}$ will be a coaction on a semicategory $\mathcal{H}$, and $\mathcal{C}$ should correspond to an ideal of $\mathcal{H}$. We note that reversing this situation, that is, beginning with a $H$-comodule semicategory $\mathcal{H}$ which contains an ideal $\mathcal{C}$ which is a category, we always get an induced coaction. 

\begin{proposition}Let $\mathcal{H}$ be an $H$-comodule semicategory with coaction $\gamma$. If $\mathcal{C}$ is a subcategory that is an ideal of $\mathcal{H}$, then $\mathcal{C}$ is a symmetrical partial $H$-comodule category with partial coaction $\rho$ given by 
\begin{eqnarray*}
\rho_{(x,y)}(f)&=&(1_{_y\mathcal{C}_y}\otimes 1_H)\gamma_{(x,y)}(f)\\
&=&\gamma_{(x,y)}(f)(1_{_x\mathcal{C}_x}\otimes 1_H),
\end{eqnarray*}  
for every $f\in\,_y\mathcal{C}_x$.
\end{proposition}

The identities $1_{_x\mathcal{C}_x}$ play a distinguished role in the construction of the globalization. 
This family defines a central idempotent of $\mathcal{H}$.

\begin{definition}[\cite{alves2}]
A \emph{central idempotent} in a linear semicategory $\mC$ is an idempotent natural transformation $e$ of the identity functor $Id_{\mC}$ to itself. That is, it is a family $e = \{_{x}e_{x}\} \in {}_x\mC_x$, indexed by $x\in \mC_{0}$, of idempotents of the endomorphism algebras ${}_x\mC_x$ such that, for every ${}_yf_x \in {}_y\mC_x$,
\begin{equation}\label{central}
{}_ye_y \circ {}_yf_x = {}_yf_x \circ {}_xe_x.
\end{equation}
\end{definition}

We introduce now what the concept of a globalization of a partial coaction on a linear category. The next definition combines notions from \cite{alves2} and \cite{alves3}.

\begin{definition}\label{def.glob.category} Let $H$ be a Hopf algebra and $\mathcal{C}$ a partial $H$-comodule category with symmetrical partial coaction $\rho$. A globalization, or enveloping coaction, of $\rho$ is a pair $(\mathcal{B},F)$, where
\begin{enumerate}[\normalfont(i)]
	\item $\mathcal{B}$ is an $H$-comodule semicategory, with coaction $\gamma$;
	\item $F: \mathcal{C}\to\mathcal{B}$ is a faithful $\mathcal{C}_0$-semifunctor and $F(\mathcal{C})$ is the ideal of $\mathcal{B}$ generated by the central idempotent $e=\{F(\,_x1_x)\}_{x\in\mathcal{C}_0}$;
	\item $\mathcal{B}$ is generated by $F(\mathcal{C})$ as $H$-comodule semicategory: for every $x,y\in \mathcal{C}_0=\mathcal{B}_0$, the morphism space $_y\mathcal{B}_x$ is generated as vector space by monomials $w^k\cdots w^1$, where for each $i$ we have ${w^i=h^*_i\rightharpoonup f^i=(I\otimes h^*)\gamma(f^i)}$, with $h^*_i\in H^*$, $f^i\in\,_{x^{i+1}}\mathcal{C}_{x^i}$ and $x^1=x$, $x^k=y$.
	\item $F$ intertwines the partial coaction on $\mathcal{C}$ and the induced partial coaction on $F(\mathcal{C})$, i.e., for every $f\in \,_y\mathcal{C}_x$, we have that 
    \[
    (F\otimes I)\rho_{(x,y)}(f)=(F(\,_y1_y)\otimes 1_H)\gamma_{(x,y)}(F(f))=\gamma_{(x,y)}(F(f))(F(\,_x1_x)\otimes 1_H).
    \]
\end{enumerate}
\end{definition}
Let $\mathcal{C}$ be a partial $H$-comodule category with symmetrical partial coaction $\rho$.
Consider the category $\mathcal{C}\otimes H$ where $_y(\mathcal{C}\otimes H)_x=\,_y\mathcal{C}_x\otimes H$, and composition is defined tensor-wise as 
\[
({}_zg_y \otimes h)({}_yf_x \otimes k)
= {}_zg_y \ {}_yf_x \otimes hk.
\]
Note that $\mathcal{C}\otimes H$ is a $H$-comodule category with coaction $\delta_{(x,y)}=I_{_y\mathcal{C}_x}\otimes \Delta$, for every $x,y\in\mathcal{C}_0$. Let $\mathcal{B}$ be the subsemicategory generated by $\rho(\mathcal{C})$, where $\rho$ is seen as a $\mathcal{C}_0$-semifunctor. In other words, for every $x,y \in\mathcal{C}_0$, the set $_y\mathcal{B}_x$ is generated as vector space by monomials $w^k\cdots w^1$, where $w^i=h^*_i\rightharpoonup f^i$, with $h^*_i\in H^*$, $f^i\in\,_{x^{i+1}}\mathcal{C}_{x^i}$ and $x^1=x$, $x^k=y$. Note that, in this case, the action $\rightharpoonup$ is given by 
\begin{eqnarray*}
	h^*\rightharpoonup \rho_{(x,y)}(f)&=&(I_{_y\mathcal{C}_x}\otimes h^*)\delta_{x,y}(\rho_{(x,y)}(f))\\
	&=&(I_{_y\mathcal{C}_x}\otimes h^*)(I_{_y\mathcal{C}_x}\otimes \Delta)\rho_{(x,y)}(f)\\
	&=&\sum f^{(0)}\otimes (f^{(1)})_{(1)}h^*((f^{(1)})_{(2)}).
\end{eqnarray*}

\begin{theorem}\label{globalization.category}
Let $\mathcal{C}$ be a partial $H$-comodule category with symmetrical partial coaction $\rho$. The pair $(\mathcal{B},\rho)$, as above, is a globalization of the symmetrical partial coaction $\rho$.
\end{theorem}
\begin{proof} Since the other axioms of globalization hold, by hypothesis, and $\rho$ is a faithful functor because of the property (i) of partial coaction (Definition \ref{def:partial.coaction}), we only need to verify that $\rho(\mathcal{C})$ is the ideal of $\mathcal{B}$ generated by the central idempotent $e=\{\rho_{(x,x)}(\,_x1_x)\}_{x\in\mathcal{C}_0}$. In fact, let $w=w^k\cdots w^1$ be a monomial in $_y\mathcal{B}_x$. First, we must prove that $\rho_{(y,y)}(\,_y1_y)w$ and $w\rho_{(x,x)}(\,_x1_x)$ lie in $\rho_{(x,y)}(\,_y\mathcal{C}_x)$. In fact, we will prove that $\rho_{(y,z)}(f)w\in \rho_{(x,z)}(\,_z\mathcal{C}_x)$ for every $f\in \,_z\mathcal{C}_y$ and every $z\in \mathcal{C}_0$. Let $h^*\in H^*$, $g\in\,_y\mathcal{C}_x$ and $f\in\,_z\mathcal{C}_y$, then
\begin{eqnarray*}
	\rho_{(y,z)}(f)(h^*\rightharpoonup g)&=&\sum f^{(0)}g^{(0)}\otimes f^{(1)}(g^{(1)})_{(1)}h^*((g^{(1)})_{(2)})\\
	&=&\sum f^{(0)}\,_y1_y^{(0)}g^{(0)}\otimes f^{(1)}\,_y1_y^{(1)}(g^{(1)})_{(1)}h^*((g^{(1)})_{(2)})\\
	&=& \sum f^{(0)}g^{(0)(0)}\otimes f^{(1)}g^{(0)(1)}h^*(g^{(1)})\\
	&=&\sum (fg^{(0)})^{(0)}\otimes (fg^{(0)})^{(1)}h^*(g^{(1)})\\
	&=&\sum \rho_{(x,z)}(fg^{(0)}h^*(g^{(1)}))\,\in\rho_{(x,z)}(_z\mathcal{C}_x).
\end{eqnarray*}
The fact that $w\rho_{(z,x)}(f)\in \rho_{(z,y)}(\,_y\mathcal{C}_z)$ for every $f\in \,_x\mathcal{C}_z$ and every $z\in \mathcal{C}_0$, is proved analogously using the fact that $\rho$ is symmetrical. Hence, particularly, $\rho_{(y,y)}(\,_y1_y)(h^*\rightharpoonup g)=\sum \rho_{(x,y)}(g^{(0)}h^*(g^{(1)}))=(h^*\rightharpoonup g)\rho_{(x,x)}(\,_x1_x)$, which means that $\rho(\mathcal{C})$ is the ideal of $\mathcal{B}$ generated by the central idempotent $e=\{\rho_{(x,x)}(\,_x1_x)\}_{x\in\mathcal{C}_0}$.
\end{proof}

%%%%%%%%%%%%%%

\subsection{Globalization of a partial coaction on an algebra with local units}

Here we discuss the globalization of a partial coaction for algebras with local units. Of course, Theorem \ref{theorem.globalization} already guarantees the existence of a globalization for a given partial coaction; the difference in the current context is that now we can show that any partial coaction can be identified with the restriction of a global coaction to an ideal.

\begin{definition}[\cite{ferrazza1}]
	Let $A$ be an algebra with local units and $S$ be an s.l.u. of $A$. A subfamily $T\subseteq S$ which is also an s.l.u. of $A$, will be called a subsystem of local units of $S$.
\end{definition}

Consider an $H$-comodule algebra $B$ with $H$-coaction $\gamma : B \to B \otimes H$ and an ideal $A$ of $B$. Ideally, if $A$ is an algebra with local units, we expect that the global $H$-coaction might restrict to a partial $H$-coaction on $A$. 

In \cite{alves3}, considering $A$ to be a \emph{unital} ideal of an algebra $B$, the induced partial coaction is
\[
\rho(a) = (1_A \otimes 1_H) \gamma (a).
\]
A reasonable adaptation for an algebra with local units would be to define
\[
\rho(a) = (e \otimes 1_H) \gamma (a)
\]
where $e$ is a local unit such that $ea = a$. Given that $e$  is not unique and depends on $a$, a potential problem here is that $\rho$ might not even be well-defined.   

We now introduce a compatibility condition between the $H$-coaction on $B$ and the system of local units of $A$. This condition is sufficient to define an induced partial coaction, and it turns out that every partial coaction is isomorphic to an induced partial coaction (Proposition \ref{prop:glob.local.units} and Theorem \ref{thm:global.local.units}).

\begin{definition} \label{def:compatibility}
        Let $B$ be an $H$-comodule algebra with coaction $\gamma:B\to B\otimes H$ and $A$ be an ideal with local units of $B$. We will say that $\gamma$ is compatible with $A$ if for every $a\in A$ there exist subsystems of local units $L_a, R_a$ such that $(e\otimes 1_H)\gamma(a)=\gamma(a)(f\otimes 1_H)$, for all $e\in L_a, f\in R_a$.
\end{definition}
\begin{proposition}\label{coacparcialideal}
    Let $B$ be an $H$-comodule algebra with coaction $\gamma:B\to B\otimes H$ and $A$ be an ideal with local units of $B$.
    If $\gamma$ is compatible with $A$ then the mapping $\rho(a)=(e\otimes 1_H)\gamma(a)$, with $e\in L_a$, defines a partial coaction $\rho:A\to A\otimes H$.
\end{proposition}

\begin{proof}
    Following similar calculations presented in the proof of Proposition 50 of \cite{ferrazza1}, we have that:
    \begin{enumerate}[\normalfont(i)]
        \item For every $e,e'\in L_a$ such that $ee'=e'e$, we can add $ee'$ to $L_a$;
        \item If $e\in L_a$ and $x$ is a local unit such that $x\geq e$, then we can add $x$ to $L_a$;
        \item If $a_1,\dots,a_n\in A$, we can choose $L_{a_i}=L$ for all $i=1,\dots,n$. Analogously, we can choose $R_{a_i}=R$ for all $i=1,\dots,n$.
    \end{enumerate}
    Clearly $\rho$ is well defined, since $(e\otimes 1_H)\gamma(a)=(e'\otimes 1_H)\gamma(a)$ for every $e,e'\in L_a$. Also, for every $e\in L_a$, we have
    \[
    (I\otimes \varepsilon)\rho(a)=(I\otimes \varepsilon)(e\otimes 1_H)\gamma(a)=ea.
    \]
    Since $L_a$ is a subsystem of local units, there is an $x\in L_a$ such that $xa=a$; therefore we must have $ea=a$ for every $e\in L_a$.

    For the second axiom of partial coaction, choosing $e\in L_a$ such that $be=b$, $e'$ in a subsystem of local units $L$ which works for every term ``$ba^{(0)}$" and for $b$, and choosing $e''$ in $L_a$ such that $\rho(b)(e''\otimes 1_H)=\rho(b)$, we have that
    \begin{eqnarray*}
        (\rho\otimes I)((b\otimes 1_H)\rho(a))&=&\sum(\rho\otimes I)(bea^{(0)}\otimes a^{(1)})\\
        &=&\sum(\rho\otimes I)(ba^{(0)}\otimes a^{(1)})\\
        &=&\sum(e'\otimes 1_H\otimes 1_H)(\gamma\otimes I)(ba^{(0)}\otimes a^{(1)})\\
        &=&\sum((e'\otimes 1_H)\gamma(b)\otimes I)(\gamma\otimes I)(a^{(0)}\otimes a^{(1)})\\
        &=&((e'\otimes 1_H)\gamma(b)\otimes 1_H)(\gamma\otimes I)\gamma(a)\\
        &=&(\rho(b)\otimes 1_H)(e''\otimes 1_H\otimes 1_H)(I\otimes \Delta)\gamma(a)\\
        &=&(\rho(b)\otimes 1_H)(I\otimes \Delta)((e''\otimes 1_H)\gamma(a))\\
        &=&(\rho(b)\otimes 1_H)(I\otimes \Delta)\rho(a).
    \end{eqnarray*}
\end{proof}

\begin{proposition} \label{prop:glob.local.units} Let $H$ be a Hopf algebra, $A$ a partial $H$-comodule algebra with  coaction $\rho$ and $S$ an s.l.u. for $A$. The pair $(B,\theta)$ is an enveloping coaction, or globalization, for the partial coaction $\rho$ if and only if
\begin{enumerate}[\normalfont(1)]
	\item $B$ is an $H$-comodule algebra with coaction $\gamma$;
	\item $\theta: A\longrightarrow B$ is a monomorphism of algebras;
	\item $\theta(A)$ is an ideal of $B$;
	\item $(\theta(e)\otimes 1_H)\gamma(\theta(a))=\gamma(\theta(a))(\theta(f)\otimes 1_H)$ for every pair of local units $e, f\in S$ such that $(e \otimes 1_H)\rho(a)=\rho(a)=\rho(a)(f\otimes 1_H)$, and
	\begin{eqnarray*}
	(\theta\otimes I)\rho(a)=(\theta(e)\otimes 1_H)\gamma(\theta(a)),
	\end{eqnarray*}
	for every $e\in S$ such that $(e \otimes 1_H)\rho(a)=\rho(a)$;
	\item $B$ is generated by $\theta(A)$ as $H$-comodule.
\end{enumerate}
Item $(4)$ means that $\theta$ is a morphism of partial coactions.
\end{proposition}
\begin{proof}
    Note that items (i), (ii), (iii), and (v) of Definition \ref{def.globalization} are the same as items (1), (2), (3), and (5), respectively. Now, if item (iv) of Definition \ref{def.globalization} holds, then we may choose $b=e$ such that $(e \otimes 1_H)\rho(a)=\rho(a)$, which gives
    \[
    (\theta\otimes I)\rho(a)=(\theta(e)\otimes 1_H)\gamma(\theta(a)).
    \]
    Also, we may choose $b=f$ such that $\rho(a)=\rho(a)(f\otimes 1_H)$, which gives
    \[
    (\theta\otimes I)\rho(a)=\gamma(\theta(a))(\theta(f)\otimes 1_H).
    \]
    Hence item (4) holds true.
    Conversely, if (4) holds, then choosing $e\in S$ such that $(e \otimes 1_H)\rho(a)=\rho(a)$ and $be=b$, we have that 
    \begin{eqnarray*}
        (\theta\otimes I)((b\otimes 1_H)\rho(a))&=&(\theta(b)\otimes 1_H)(\theta\otimes I)\rho(a)\\
        &=&(\theta(b)\otimes 1_H)(\theta(e)\otimes 1_H)\gamma(\theta(a))\\
        &=&(\theta(b)\otimes 1_H)\gamma(\theta(a))
    \end{eqnarray*}
    The equality $(\theta\otimes I)(\rho(a)(b\otimes 1_H))=\gamma(\theta(a))(\theta(b)\otimes 1_H)$ is proved analogously.
\end{proof}

\begin{remark}
Note that if $(B,\theta)$ is a globalization for a partial coaction on an algebra with local units, then $(B,\theta)$ is a globalization in the general case. Also, if $(B,\theta)$ is a globalization for a partial coaction $\rho:A\to A\otimes H$ and $A$ has local units, then for every $a\in A$ the sets
\[
L_{\theta(a)}=\{\theta(e)\,|\,(e\otimes 1_H)\rho(a)=\rho(a)\}\quad\text{and}\quad
R_{\theta(a)}=\{\theta(f)\,|\,\rho(a)(f\otimes 1_H)=\rho(a)\}
\]
satisfy the compatibility condition of Definition \ref{def:compatibility} for the ideal $\theta(A)$ of $B$. In this case, we see explicitly how the partial coaction on $\theta(A)$ is associated with the coaction on $B$: for each $a \in A$ we have
\[
\Phi(\theta(a))=(\theta(e)\otimes 1_H)\gamma(\theta(a))
\]
for some local unit $e\in S$; hence we may say that the partial coaction $\Phi$ is induced by the coaction $\gamma$.
\end{remark}

\begin{theorem} \label{thm:global.local.units}
Let $H$ be a Hopf algebra, $A$ a partial $H$-comodule algebra with local units with symmetrical partial coaction $\rho$, and let $S$ be an s.l.u. for $A$. Then $(B,\rho)$ is an enveloping coaction for $\rho$, where $B$ is the subcomodule algebra of $A\otimes H$ generated by $\rho(A)$.
\end{theorem} 
\begin{proof}
This is a consequence of Theorem \ref{theorem.globalization}.
\end{proof}

We conclude this work by showing how to relate the globalization of a $S$-categorizable partial coaction on $A$ with the globalization of $\mC^S(A)$ and, conversely, the globalization of a partial $H$-comodule category $\mC$ with the globalization of the algebra $a(\mC)$.

\begin{theorem} \label{thm:global.S-categorical.to.category}
Let $A$ be an algebra with s.l.u. $S$ and $H$ a finite-dimensional Hopf algebra. If $\rho:A\to A\otimes H$ is an $S$-categorizable symmetrical partial coaction and $(B,\theta)$ is a globalization with coaction $\gamma$, then the induced partial $H$-coaction on $\mathcal{C}^S(A)$ has a globalization $(\mathcal{B},F)$, where:
\begin{enumerate}[\normalfont(1)]
    \item $\mathcal{B}_0=\mathcal{C}^S(A)_0$, and $\mathcal{B}(e,f)$ is generated as a vector space by monomials $w^k\cdots w^1$, where $w^i=h_i^*\rightharpoonup a^i=(I\otimes h_i^*)\gamma(a^i)$, with $h_i^*\in H^*$, $a^i\in\mathcal{B}(e_i,e_{i+1})$, and $e_1=e$, $e_k=f$ (note that $\mathcal{B}(e,f)\subset B$);
    \item the coaction on the semicategory $\mathcal{B}$ is given by $\{\gamma(e,f)=\gamma|_{\mathcal{B}(e,f)}\}$;
    \item $F:\mathcal{C}^S(A)\to \mathcal{B}$ is the semifunctor given by $F_0(e)=e$ and $F_1(fae)=\theta(fae)$.
\end{enumerate}
\end{theorem}

\begin{proof}
Axiom (i) of Definition \ref{def.glob.category} follows from $B$ being an $H$-comodule algebra with coaction $\gamma$. Axiom (ii) follows from the axioms of Definition \ref{def.globalization} and similar calculations as in the proof of Theorem \ref{globalization.category}. Axiom (iii) follows directly from the construction of each $\mathcal{B}(e,f)$. Finally, axiom (iv) of Definition \ref{def.glob.category} follows from item (4) of Definition \ref{def.globalization} and from the fact that the partial coaction $\rho$ is $S$-categorizable.
\end{proof}

\begin{theorem} \label{thm:global.category.to.S-categorical.}
Let $\mathcal{C}$ be a partial $H$-comodule category with coaction $\{\rho(x,y):\,_y\mathcal{C}_x\to\,_y\mathcal{C}_x\otimes H\}_{x,y\in\mathcal{C}_0}$ and let $(\mathcal{B}, F)$ be a globalization for it. Then $(a(\mathcal{B}),\theta)$ is a globalization for the induced partial coaction on $a(\mathcal{C})$, where $\theta: a(\mathcal{C})\to a(\mathcal{B})$ is given by
\[
\theta\bigl((\,_yf_x)_{x,y\in\mathcal{C}_0}\bigr)=(F(\,_yf_x))_{x,y\in\mathcal{B}_0}.
\]
\end{theorem}

\begin{proof}
For item (i) of Definition \ref{def.globalization}, it is enough to consider $\gamma: a(\mathcal{B})\to a(\mathcal{B})\otimes H$ defined by $\gamma((\,_yf_x)_{x,y\in\mathcal{B}_0})=(\gamma(x,y)(\,_yf_x))_{x,y\in\mathcal{B}_0}$. Item (ii) is a consequence of $F$ being a faithful $\mathcal{C}_0$-semifunctor. Items (iii), (iv), and (v) of Definition \ref{def.globalization} are consequences of items (ii), (iv), and (iii) of Definition \ref{def.glob.category}, respectively.
\end{proof}

\bibliographystyle{plain}
\bibliography{refs}

\end{document}